\documentclass[12pt]{article}
\usepackage{WideMargins2}

\def\GCD{{{\rm GCD}}}
\def\Norm{{{\rm N}}}

\def\F{{{\mathbb F}}}
\def\Z{{{\mathbb Z}}}

\def\cj#1{{{\overline{#1}}}}

\def\deg{{{\rm deg}}}
\def\ang#1{{{\langle #1 \rangle}}} 
\def\del#1{{{\ang{\mu_{#1}}}}}
\def\delstar#1{{{\ang{\mu_{#1}^*}}}}

\let\varep=\varepsilon
\let\x=\times
\let\s=\sigma

\let\ee = e
\def\ie{{ \it i.e.}}
\def\ang#1{{{\langle #1 \rangle}}}
\def\set#1{{\left\{\,#1\,\right\}}}
\def\jacobi#1#2{{{\left(\frac{#1}{#2}\right)}}}
\def\choose#1#2{{{\left({#1}\atop {#2}\right)}}}
\def\calT{{{\cal T}}}
\def\calZ{{{\cal Z}}}
\def\calA{{{\cal A}}}
\def\calB{{{\cal B}}}

\def\calS{{{\cal S}}}
\def\calN{{{\cal N}}}
\def\typ{{{\rm Type}}}

\title{\Large\bfseries  Permutation properties of
Dickson and Chebyshev polynomials with connections to number theory}
\author{Antonia W.~Bluher\\ National Security Agency \\ tbluher@access4less.net}
\date{\today}
\begin {document}
\fancytitle

\begin{abstract}
The $k$th Dickson polynomial of the first kind, $D_k(x)\in\Z[x]$, is determined by the formula: $D_k(u+1/u) = u^k + 1/u^k$,
where $k\ge 0$ and $u$ is an indeterminate.  
These polynomials are closely related to Chebyshev polynomials and 
have been widely studied.
Leonard Eugene Dickson proved in 1896
that $D_k(x)$ is a permutation polynomial on $\F_{p^n}$, $p$ prime, 
if and only if $\GCD(k,p^{2n}-1)=1$, and his result easily carries over to Chebyshev polynomials
when $p$ is odd.
This article continues on this theme, as we find special subsets of
$\F_{p^n}$ that are stabilized or permuted by Dickson or Chebyshev polynomials.
Our analysis also leads to a
factorization formula for Dickson and Chebyshev polynomials and some
new results in elementary number theory.
For example, we show that if $q$ is an odd prime power, then 
$\prod \{ a \in \F_q^\x : \text{$a$ and $4-a$ are nonsquares} \} = 2$.
\end{abstract}

\section{Introduction} \label{sec:intro}

This article finds certain natural subsets of $\F_q^\x$ 
that are stabilized or permuted by Dickson or Chebyshev polynomials, where $q$ is a prime power.
This analysis leads to discovery of 
a factorization formula for Dickson and Chebyshev polynomials and some 
new results in elementary number theory.  

The $k$th Dickson polynomial of the first kind, $D_k(x)\in\Z[x]$, is determined by the recursion
$D_0(x)=2$, $D_1(x)=x$, and $D_{k+2}(x)=xD_{k+1}(x)-D_k(x)$ for $k\ge0$. It can be shown by induction
that\footnote{There is a more general definition of Dickson polynomial of the first kind that depends on a parameter~$a$.
In the more general definition, $D_k(x,a)$ is determined
by the formula $D_k(u+a/u,a) = u^k + (a/u)^k$.  However, for this article we are interested only in the case $a=1$.}
$$D_k(u+1/u) = u^k + 1/u^k,$$
where $k\ge 0$ and $u$ is an indeterminate. In fact, this functional equation determines $D_k$ uniquely and hence
can serve as an alternate definition.  Substituting $-u$ for $u$ in the functional equation implies the well-known fact that $D_k(-x)=(-1)^k D_k(x)$.

Any polynomial in $\Z[x]$ determines a function from $\F_q$ to itself.
Leonard Eugene Dickson proved in his Ph.D.\ thesis (1896)
that $D_k(x)$ permutes $\F_{q}$ if and only if $\GCD(k,q^2-1)=1$.
The Chebyshev polynomial of the first kind, $T_k(x)$, is related to the Dickson polynomial
by the formula $2T_k(x) = D_k(2x)$.  (Here $T$ stands for Tchebycheff, which is the French transliteration of Chebyshev.)
If 2 is invertible (e.g. in odd characteristic fields), $D_k$ and $T_k$ are essentially equivalent (in fact, conjugate), and in particular
Dickson's result carries over to Chebyshev polynomials when $q$ is odd.
However, in characteristic~2, the Dickson polynomials have a theory while the Chebyshev polynomials do not.
For this reason, we focus mainly on Dickson polynomials, although we do summarize our results in terms of Chebyshev polynomials
in Section~\ref{sec:Chebyshev}.

If $q=2^n$, then it was shown by Dillon and Dobbertin \cite{DD} that $D_k(\calT_0)\subset \calT_0\sqcup \{0\}$
and $D_k(\calT_1) \subset \calT_1 \sqcup \{0\}$, where $\calT_j= \left\{ a \in \F_q^\x : \Tr(1/a)=j \right\}$.
A main theme of this article is to produce odd-characteristic analogues of this result, as described 
below.  

If $q$ is odd and $a\in\F_q$ then let $\jacobi aq$ denote the Legendre symbol:
\begin{equation}
 \jacobi a q = \begin{cases} 1 & \text{if $a$ is a square in $\F_q^\x$} \\
$-1$ & \text{if $a$ is a nonsquare in $\F_q^\x$} \\ 0 & \text{if $a=0$}. \end{cases}  \label{jacobiDef}
\end{equation}
For $\varepsilon_1,\varepsilon_2 \in \{1,-1\}$ and $\lambda \in \F_q^\x$, define
\begin{equation}\calA_\lambda^{\varepsilon_1,\varepsilon_2}
= \left\{\, u \in \F_q : \jacobi {u-\lambda}q = \varepsilon_1,\quad
\jacobi {u+\lambda}q = \varepsilon_2 \,\right\}\label{Aij}. \end{equation}
We often write $+,-$ instead of 1, $-1$, {\it e.g.} $\calA^{+-}_\lambda$ instead of $\calA^{1,-1}_\lambda$.
Note that  $\F_q$ is a disjoint union:
\begin{equation*} \F_q = \calA_\lambda^{++} \sqcup \calA_\lambda^{+-} \sqcup \calA_\lambda^{-+} \sqcup \calA_\lambda^{--} \sqcup \{\lambda,-\lambda\}. \end{equation*}
We prove that $D_k(x)$ permutes each set $\calA_2^{\pm,\pm}$ and
$T_k(x)$ permutes each set $\calA_1^{\pm,\pm}$ when $\GCD(k,q^2-1)=1$. 
Even when $k$ is not assumed to be relatively
prime to $q^2-1$, we can still
say a lot about the image sets $D_k(\calA_2^{\varepsilon_1,\varepsilon_2})$
or $T_k(\calA_1^{\varep_1,\varep_2})$. The following theorem about Dickson polynomials will be
proved in Section~\ref{sec:refinements}. A corresponding theorem for Chebyshev polynomials is
stated in Theorem~\ref{thm:ChebyshevRefinements}.

\begin{theorem}  \label{thm:Aij} ($q$ odd.) Let $\varepsilon=(-1)^{(q-1)/2}$, so $q \equiv \varepsilon \pmod 4$ and $\varepsilon = \jacobi{-1}q$.
For any $k\ge1$, $D_k(2)=2$, $D_k(-2)=(-1)^k \cdot 2$, and
\begin{enumerate}
\item[{\it (i)}] $D_k(\calA_2^{++}) \subset \calA_2^{++} \sqcup \{2,-2\varepsilon \}$ and $D_k(\calA_2^{-+}) \subset \calA_2^{-+} \sqcup \{2,2\varepsilon\}$.  
\item[{\it (ii)}] If $k$ is odd then $D_k(\calA_2^{\varepsilon,-})\subset \calA_2^{\varepsilon,-} \sqcup \{-2\}$ and
$D_k(\calA_2^{-\varepsilon,-}) \subset \calA_2^{-\varepsilon,-}$. 
\item[{\it (iii)}] If $k$ is even then $D_k(\calA_2^{--}) \subset \calA_2^{++}\sqcup \{2,-2\varepsilon \}$ and 
$D_k(\calA_2^{+-})\subset \calA_2^{-+}\sqcup \{2,2\varepsilon \}$. 
\item[{\it (iv)}] If $(q-1)/2$ divides $k$, then $D_k(\F_q) \subset \calA_2^{+-} \sqcup \calA_2^{-+} \sqcup \{2,-2\}$.
If in addition $k$ is even, then $D_k(\F_q) \subset \calA_2^{-+} \sqcup \{2,-2\}$.
\item[{\it (v)}] If $(q+1)/2$ divides $k$, then $D_k(\F_q) \subset \calA_2^{++} \sqcup \calA_2^{--} \sqcup \{2,-2\}$. 
If in addition $k$ is even, then $D_k(\F_q) \subset \calA_2^{++} \sqcup \{2,-2\}$.
\item[{\it (vi)}] If $\calA_2^{\varepsilon_1,\varepsilon_2}$ is nonempty, then 
$D_k$ permutes $\calA_2^{\varepsilon_1,\varepsilon_2}$ iff $\GCD(k,d^{\varepsilon_1,\varepsilon_2})=1$, where
$d^{++}=(q-1)/2$, $d^{+-}=q+1$, $d^{-+}=(q+1)/2$, $d^{--} = q-1$. 
\end{enumerate}
\end{theorem}

Another result of this article concerns the factorization of Dickson polynomials. Assuming $q$ is odd,
let $\varepsilon=(-1)^{(q-1)/2}$ and $m=(q-\varepsilon)/4$.  Then $\jacobi 2q = (-1)^m$ (see \cite[eq. (26)]{wilson-like}),
so $D_m(-2)=\jacobi2q\cdot 2$.
We prove in Theorem~\ref{thm:fijDE} that
$$ D_m(x) = \prod_{b \in \calA_2^{-\varepsilon,-}} (x-b). $$
As shown in Theorem~\ref{thm:wilson-like}, evaluating the above formula at $x=2$ and $x=-2$ gives
\begin{equation}
\prod \left\{ a \in \F_q^\x : \text{$a$ and $4-a$ are nonsquares} \right\} = 2,
\label{T40}
\end{equation}
\begin{equation}
\prod \left\{ a \in \F_q^\x : \text{$-a$ and $4+a$ are nonsquares} \right\} = \jacobi 2 q\cdot  2.\label{T04}
\end{equation}
These formulas are striking for their similarity to Wilson's Theorem, which states that 
the product of all nonzero elements of a field $\F_q$ is $-1$.  
We decided to continue the development of Wilson-like theorems in a separate 
article \cite{wilson-like}, as it was too lengthy to include in the current work; see also
\cite{structure}, which is a short summary of results from both this article and \cite{wilson-like}.
In \cite{wilson-like}, explicit formulas are found for
$\prod \calT_{j,\ell}^{\varep_1,\varep_2}$
for all $j,\ell\in\F_q$ and for all $\varepsilon_1,\varepsilon_2\in \{1,-1\}$,
where
\begin{equation}
\calT_{j,\ell}^{\varep_1,\varep_2}=
\left\{ a \in \F_q^\x : \jacobi{j-a}q=\varepsilon_1,\ \jacobi{\ell+a}q=\varepsilon_2\right\}.\label{TjellDef}
\end{equation}
Formulas (\ref{T40}), (\ref{T04}) are the special cases: $\prod\calT_{4,0}^{--}=2$
and $\prod\calT_{0,4}^{--}=\jacobi2q \cdot 2$.

Several other new results of a number-theoretic flavor are given in 
Sections~\ref{sec:numberTheory} and~\ref{sec:sigmak}.
For example, Proposition~\ref{prop:D2}({\it i}) states that if $q$ is
odd then
$b \mapsto b^2-2$ is a permutation of $\left\{b\in\F_q : \jacobi{2-b}q=-1,\ \jacobi{2+b}q = 1 \,\right\}$, and the inverse of this permutation is
$$b \mapsto \prod \{(b-a) : \text{$2-a$ and $2+a$ are nonsquares} \}.$$
Another new formula is the identity for all $c \in\F_q$ ($q$ odd):
\begin{eqnarray}\label{jacc}
&&\prod\left\{c-a : a \in\F_q,\ \jacobi{a(a+4)}q= 1 \right\} + \\ \nonumber
&&\prod\left\{c-b : b \in\F_q,\ \jacobi{b(b+4)}q= -1 \right\} = \jacobi{c}q. 
\end{eqnarray}
Proposition~\ref{prop:sqrtc} demonstrates that if $q=2^n$, then for any $c \in \F_q$,
\begin{equation} \label{sqrtc}
 \prod \{ c + 1/a : a \in \F_q^\x,\ \Tr(a) = 0 \} 
+ \prod \{ c + 1/b : b \in \F_q^\x,\ \Tr(b)=1 \} = c^{1/2}. 
\end{equation}
Theorem~\ref{thm:sigmak} gives closed formulas for all the
elementary symmetric functions of four special subsets of $\F_q$
when $q$ is odd, providing an additional way to generalize
Wilson-like formulas.

Formulas~(\ref{T40}) and (\ref{sqrtc}), although discovered using Dickson
polynomials, also admit direct proofs. Short proofs were communicated
to the author by
Richard Stong and John Dillon and are presented in
Section~\ref{sec:oneline}.  

\medskip
\noindent{\it Generalities. } 
Throughout this article, $p$ denotes a prime and $q=p^n$.  
Here, $p=2$ is allowed.  We use the notation:
$$\ang{u} = u + 1/u,$$
where $u\ne 0$.
Then the defining equation for $D_k$ is
$$D_k(\ang u) = \ang{u^k}.$$
In other words, there is a commutative square:
\begin{equation*}
\xymatrix{ u \ar[r]^{x^k} \ar[d] &  u^k \ar[d] \\
\ang{u} \ar[r]^{D_k} & \ang{u^k} }
\end{equation*}
Observe that $\ang{u^k}\ang{u^\ell} = \ang{u^{k+\ell}} + \ang{u^{k-\ell}}$, and so (setting
$x=\ang u$) one deduces the well-known formula if $k \ge \ell$,
\begin{equation}D_k(x) D_\ell(x) = D_{k+\ell}(x) + D_{k-\ell}(x).\label{DkDl} \end{equation}

\medskip
\noindent{\it Notation.}
We use much notation that is by now very standard: 
$\Z$ denotes the ring of integers, ${\mathbb Q}$ denotes the field of 
rational numbers, and
$\F_q$ denotes the unique field with $q$ elements.
If $R$ is a ring then $R^\x$ denotes the group of invertible
elements of $R$ and $R[x]$ is the polynomial ring over $R$; if $S$ is a set then $|S|$ denotes its cardinality. The symbol for union is $\cup$.
$S_1\sqcup S_2 \sqcup \cdots \sqcup S_k$ denotes 
the union of sets $S_1,\ldots,S_k$ that are known to be pairwise disjoint. 
If $F$ is a field, then $\cj F$ denotes its algebraic closure. 
We sometimes write iff or $\iff$ as an abbreviation for ``if and only if''.
We assume the reader has background in finite fields, as can be found for example in \cite[Chapter 2]{LN}. Throughout, $q=p^n$ denotes a prime power.
If $q$ is odd, then we set $\varepsilon=(-1)^{(q-1)/2}=\jacobi{-1}q$ and $m=(q-\varep)/4\in\Z$.

\section{Analysis of Dickson polynomials over finite fields} 
\label{sec:analysis}

For $d\in\Z$ not divisible by $p$, define (in characteristic~$p$)
\begin{equation} \mu_d = \left\{a \in \cj\F_p^\x : a^d = 1 \right\},\qquad 
\del d = \left\{\ang{a} : a \in \mu_d \right\}. \label{mu_delta_def}
\end{equation}
The characteristic~$p$ is not included in the notation for $\mu_d$ and $\del d$, but it will be clear from the context.
It is well known from the theory of finite fields that $\mu_d$ is a cyclic group of order~$d$ and that
$\mu_{p^n-1} = \F_{p^n}^\x$. We have $\mu_k \cap \mu_d = \mu_g$, where $g = \GCD(k,d)$.
If $uv\ne 0$, then $\ang u = \ang v \iff u^2v + v = uv^2 + u \iff (uv-1)(u-v)=0$, and so
\begin{equation} \ang u = \ang v \iff v \in \{u,1/u\}.\label{angEquality} \end{equation}
In particular, $u$ and $v$ have the same multiplicative order. From this, we see that 
\begin{equation} \del k \cap \del d = \del g,\qquad \text{where $g=\GCD(k,d)$} \label{delg}
\end{equation}
and
\begin{equation}
\ang{\mu_k} \subset \ang{\mu_d} \iff k|d,\qquad
\ang{\mu_k} = \ang{\mu_d} \iff k = d. \label{delEquality}
\end{equation}

It has been noted in \cite{Brewer}, \cite[Lemma~3.5]{Lidl}, and \cite[p.~356]{DD} that 
$u+1/u \in \F_q$ iff $u^{q+1}=1$ or $u^{q-1}=1$. In our notation, 
this observation takes the form of Lemma~\ref{lem:FqDecomposition}, below.

\begin{lemma} ($q$ even or odd.) \label{lem:FqDecomposition} $\F_q = \del{q-1} \cup \del{q+1}$. Also, $\del{q-1} \cap \del{q+1} = \{2,-2\}$. \end{lemma}

\begin{proof} 
Given $a \in \F_q$, let $u \in \cj \F_q^\x$ such that $u^2-au+1=0$. Then $u\ne 0$ and $a=u+1/u$.
By~(\ref{angEquality}),
$$a\in\F_q\iff \ang u = \ang {u^q} \iff u^q \in \{u,1/u\} \iff u \in \mu_{q-1}\cup \mu_{q+1}.$$  
If $q$ is odd, then $\del{q-1}\cap\del{q+1}=\del 2 = \{\ang 1,\ang{-1}\}=\{2,-2\}$ by (\ref{delg}).
If $q$ is even, then $\{2,-2\}=\{0\}$ and also $\del{q-1}\cap\del{q+1}=\del 1 = \{0\}$. 
\end{proof}

Define the following subsets of $\cj\F_p$:
\begin{equation} \calZ = \{2,-2\},\qquad \delstar d = \del d \setminus \calZ, \label{deltaStar} \end{equation}
where $p \nmid d$.  
Note that $\calZ = \del 2$ if $p$ is odd, and $\calZ = \{0 \} = \del 1$ if $p=2$.
By Lemma~\ref{lem:FqDecomposition}, $\F_q$ is a disjoint union:
\begin{equation} \F_q = \delstar {q-1} \sqcup \delstar {q+1} \sqcup \calZ. \label{FqDisjoint} \end{equation}
By (\ref{angEquality}), $u \mapsto \ang u $ is 2-to-1 from $\mu_d$ to $\del d$ unless $u \in \{1,-1\}$, in which case it is 1-to-1.
Noting that $\mu_d\cap\{\,1,-1\} = \{1\}$ if $d$ is odd and $\{1,-1\}$ if $d$ is even, it follows that
\begin{equation}|\delstar d|= \begin{cases} (d-1)/2 &\text{if $d$ is odd} \\  (d-2)/2 & \text{if $d$ is even} \end{cases} = \lfloor (d-1)/2 \rfloor.
\label{deltaCard} \end{equation}
In particular,
\begin{equation}
|\del {q-1}^*| = \begin{cases} q/2 - 1 & \text{if $q$ is even,} \\
(q-3)/2 & \text{if $q$ is odd,} \end{cases} \qquad
|\del {q+1}^*| = \begin{cases} q/2  & \text{if $q$ is even,} \\
(q-1)/2 & \text{if $q$ is odd.} \end{cases} \label{deltaCount}
\end{equation}
The next lemmas are well known.  As usual, assume $p\nmid d$.

\begin{lemma} \label{powerMap} Let $k\ge 1 $ and $g=\GCD(d,k)$. The map $a \mapsto a^k$ gives a $g$-to-1 map from $\mu_d$ onto $\mu_{d/g}$.
\end{lemma}
\begin{proof} 
Write $k=gk_0$ and $d=gd_0$. 
If $a \in \mu_d$, then $a^k \in \mu_{d_0}$ because
$(a^k)^{d_0} = (a^d)^{k_0}=1$.  Conversely, if $b \in \mu_{d_0}$ then
we must show $b=a^k$ for some $a\in \mu_d$. Let $c$ be a $g$-th root of $b$;
then $c \in \mu_d$. By the Euclidean 
algorithm, there are
$r,s\in\Z$ such that $rk+sd=g$. 
Then $b=c^g = c^{rk+sd}=(c^r)^k\cdot (c^d)^s = a^k$, where $a=c^r \in\mu_d$.
We have shown that $a \mapsto a^k$ is a homomorphism from $\mu_d$ onto $\mu_{d/g}$. The kernel must have order
$d/(d/g)=g$, so the map is $g$-to-1.
\end{proof}

\begin{lemma} \label{Dkdelta} Let $k\ge 1$ and $g=\GCD(d,k)$. Then $D_k(\del d) = \del {d/g}\subset \del d$. In particular, $D_k$
permutes $\del d$ if and only if $\GCD(d,k)=1$. \end{lemma}
\begin{proof} Let $\ang a \in \del d$. Then $D_k(\ang a)=\ang{a^k} \in \del {d/g}$ by Lemma~\ref{powerMap}. To prove surjectivity, let
$\ang b \in \del {d/g}$, where $b\in \mu_{d/g}$. By Lemma~\ref{powerMap}, $b=a^k$ for some $a\in\mu_d$. Then $D_k(\ang a)=\ang{a^k}=\ang b$.
We have shown that $D_k(\ang{\mu_d})=\mu_{d/g}$.  
In particular, $D_k$ permutes $\del d$ iff
$\ang{\mu_{d/g}}=\ang{\mu_d}$ iff $g=1$ by (\ref{delEquality}).
\end{proof}

\begin{proposition}   \label{prop:dicksonImage} 
Let $q=p^n$ (even or odd), $k \ge 1$,  $r = \GCD(k,q-1)$, and $s=\GCD(k,q+1)$. Then
$$D_k(\F_q) = \del {(q-1)/r} \cup \del {(q+1)/s}.$$
\end{proposition}
\begin{proof} This is immediate from Lemmas \ref{lem:FqDecomposition} 
and~\ref{Dkdelta}. 
\end{proof}

Dickson's result, that $D_k$ permutes $\F_q$ if and only if $\GCD(k,q^2-1)=1$, follows immediately from Proposition~\ref{prop:dicksonImage}, 
Lemma~\ref{lem:FqDecomposition}, and (\ref{delEquality}).

\begin{corollary} \label{cor:DkeqDl} $D_k(\F_q) = D_\ell(\F_q)$ iff $\GCD(k,q-1)=\GCD(\ell,q-1)$ and $\GCD(k,q+1) = \GCD(\ell,q+1)$.
\end{corollary}

Corollary~\ref{cor:DkeqDl} was observed by Dillon and Dobbertin \cite[pages 355--356]{DD} and was used to prove 
their beautiful theorem:

\begin{theorem} (Dillon and Dobbertin, \cite[Theorem 7]{DD}) 
Let $k=2^e+1$ or $k=2^e-1$ and $q=2^n$, and assume that $(e,n)=1$. Then $D_k(\F_q)=\F_q$ if $\GCD(k,2^{2n}-1)=1$, and otherwise $D_k(\F_q)=D_3(\F_q)$.
\end{theorem}

In summary, if $p\nmid d$, $k\ge1$, and $g = \GCD(d,k)$, then the commutative diagram shown in Figure~\ref{fig:commutativeDiagram} holds.
The map $x \mapsto x + 1/x$ is 2-to-1 unless $x=1/x$, {\it i.e.} $x\in \{1,-1\}$,
in which case $\ang x \in \calZ$.
Using this, one can work out the multiplicities of $D_k$; see \cite[Theorem 3.26]{Lidl}.
Here we just present the simplest case.

\begin{figure}[h]
\begin{equation*}
\xymatrix{ a \in \mu_d \ar[r]^{x^k} \ar[d]_{x+1/x} &  a^k \in \mu_{d/g} \ar[d]^{x+1/x}\\
\ang{a} \in \del d \ar[r]^{D_k} & \ \ang{a^k} \in \del {d/g} }
\end{equation*}
\caption{Commutative diagram}
\label{fig:commutativeDiagram}
\end{figure}

\begin{proposition} \label{Dk_multiplicity} ($q$ even or odd.) Consider $D_k : \del d \to \del {d/g}$, where $g=\GCD(d,k)$. If $y \in \delstar {d/g}=\del {d/g} \setminus \calZ$, then $y$ has
exactly $g$ preimages in $\del d$ under $D_k$, which all belong to $\delstar d$.
\end{proposition}
\begin{proof} By Lemma~\ref{Dkdelta}, the horizontal arrows in Figure~\ref{fig:commutativeDiagram} are surjective. 
Let $y \in \ang{\mu_{d/g}^*}$ and let $t$ denote the number of preimages of $y$ in $\ang{\mu_d}$; we will show that
$t=g$.
Let $\ang a\in\del d $ be a preimage of $y$, so
$y = D_k(\ang a) = \ang{a^k}$. By hypothesis, $y \not \in \calZ$. Since $D_k(\calZ) \subset \calZ$, none of the preimages of
$y$ under the map $D_k$ belong to $\calZ$, either. Note that $\ang a \not \in\calZ$ implies that $\ang a $ has exactly two preimages in $\mu_d$ under the map
$x + 1/x$ (namely $a$ and $1/a$).  Now consider the commutative diagram, 
and count the number of preimages of $y$ in the top left of the diagram in two ways.  We can pull back vertically to $\mu_{d/g}$ to obtain two preimages, and then pull back by the power map to obtain $2\x g$ preimages in $\mu_d$.
Alternatively, pull back by $D_k$ to get $t$ preimages in $\delstar d$, and then by the vertical map to obtain $2 \x t$ preimages in $\mu_d$.  Thus, $2g=2t$ and so $t=g$.
\end{proof}

\begin{corollary} \label{cor:permuteDel} ($q$ even or odd.) Suppose that $\delstar d \ne \emptyset$ (equivalently, $d>2$). Then $D_k(x)$ permutes $\delstar d$
iff GCD$(d,k)=1$.
\end{corollary}

\section{Alternative descriptions of $\delstar{q-1}$ and $\delstar{q+1}$ in even characteristic} \label{sec:evenq}

The sets $\delstar{q-1}$ and $\delstar{q+1}$ have interesting alternative descriptions. In the characteristic two case, this is well known, due to Dillon
and Dobbertin.  
Since $D_k(\del{q-1}) \subset \del{q-1}$
and $D_k(\del{q+1}) \subset \del{q+1}$, these alternative descriptions lead to interesting results 
(Proposition~\ref{DTq} and Theorem~\ref{thm:oddq})
about subsets of $\F_q$ that are preserved by Dickson polynomials. 

For $j=0,1$ and $q=2^n$ define
\begin{equation} \calT_j = \{ a \in \F_{q}^\x : \Tr_{\F_{q}/\F_2}(1/a) = j \}. \label{TqjDef}
\end{equation}
Since $\Tr_{\F_q/\F_2}$ is a surjective homomorphism, its kernel has order $q/2$, and so $|\calT_0\sqcup \{0\}|=|\calT_1|=q/2$.
$\F_q$ is a disjoint union:
\begin{equation} \F_q = \calT_{0} \sqcup \calT_{1} \sqcup \{0\},\quad \text{and $|\calT_{0}|=q/2-1$, $|\calT_{1}|=q/2$.} \label{decomp1}\end{equation}
On the other hand, by (\ref{FqDisjoint}) and (\ref{deltaCount}), $\F_q$ may also be written as a disjoint union:
\begin{equation} \F_q = \delstar{q-1} \sqcup \delstar{q+1} \sqcup \{0\},\quad 
\text{and $|\delstar{q-1}|=q/2-1$, $|\delstar{q+1}|=q/2$.}\label{decomp2}\end{equation}
The following was observed by Dillon and Dobbertin.

\begin{proposition} (\cite[page 355]{DD}) \label{calT_prop}
Let $q=2^n$. Then $\del{q-1} = \calT_{0}\sqcup\{0\}$ and $ \del{q+1}= \calT_{1}  \sqcup\{0\}$.
\end{proposition}
\begin{proof} 
Since $\calZ = \{0\}$ in characteristic~2, we see in the notation of (\ref{deltaStar}) and~(\ref{FqDisjoint}) that $\delstar{q\pm1}=\del{q\pm1} \cap \F_q^\x$,
so the proposition is equivalent to
\begin{equation}\delstar{q-1} = \calT_{0}, \qquad \delstar{q+1} = \calT_{1}  \qquad\text{when $q=2^n$.} \label{Tq01} \end{equation}
We will prove $\delstar{q-1} \subset \calT_{0}$. Since the cardinalities of these sets are the same, this will imply they are equal.
Then, (\ref{decomp1}) and (\ref{decomp2}) together will imply that $\delstar{q+1}=\calT_{1}$ as well.

Let $a \in \delstar{q-1}$;  then $a\ne 0$ and $a=r+1/r$ with $r \in \F_q^\x$.  We have $(r/a^2)(a+r+1/r)=0$, or equivalently $r/a + (r/a)^2 = 1/a^2$.
Then $\Tr_{\F_q/\F_2}(1/a)=0$, {\it i.e.} $a\in \calT_{0}$. Thus, $\delstar{q-1}\subset \calT_{0}$ as claimed.
\end{proof}

\begin{proposition} \cite[page 356]{DD}  \label{DTq} ($q$ even.)
$D_k(\calT_{0}) \subset \calT_{0} \cup \{0\}$ and $D_k(\calT_{1}) \subset \calT_{1}\cup \{0\}$.
$D_k$ permutes $\calT_{0}$ iff $\GCD(k,q-1)=1$, and it 
permutes $\calT_{1}$ iff $\GCD(k,q+1)=1$.
\end{proposition}

\begin{proof} The first sentence is immediate from 
Proposition~\ref{calT_prop} and Lemma~\ref{Dkdelta}. 
By Lemma~\ref{Dkdelta}, $D_k$ permutes $\del{q-1}$ iff $\GCD(k,q-1)=1$. 
Since $D_k(0)=D_k(\ang 1) = \ang{1^k}=0$, 
we see that $D_k$ permutes $\del{q-1}$ 
iff it permutes $\del{q-1} \setminus \{0\} = \calT_{0}$.  
A similar argument applies to $\del{q+1}$ and $\calT_{1}$.
\end{proof}

Flori and Mesnager (\cite[Section 2]{FM}) have a different proof that
$D_k(\calT_{i}) \subset \calT_{i} \cup \{0\}$.  Namely, they
compute a rational function $R_k$ such that
$1/D_k(x) + 1/x = R_k(x)^2 + R_k(x)$. From this, they deduce that if $a\in\F_q^\x$ and $D_k(a)$ is nonzero, then 
$1/D_k(a)$ and $1/a$ have the same absolute trace.
Thus, $D_k(\calT_{i})\subset \calT_{i} \cup \{0\}$. 

\begin{proposition} \label{prop:evenPermute} ($q$ even.) 
({\it i})\  
$D_k(\calT_{1})=\{0\}$ iff $q+1$ divides $k$.  
If $q>2$ then $D_k(\calT_{0})=\{0\}$ iff $q-1$ divides $k$. \\
({\it ii})\   
$D_k(b)=b^2$ for all $b \in \calT_{1}$ iff $k \equiv \pm2 \pmod {q+1}$. 
If $q>2$ then 
$D_k(a)=a^2$ for all $a \in \calT_{0}$ iff $k \equiv \pm2 \pmod {q-1}$. \\
({\it iii}) Let $e\ge 0$. If $e$ is even, then $D_{q^e-1}(a)=0$ and $D_{q^e+1}(a)=a^2$ for all $a\in\F_q$.  If $e$ is odd, then
$$D_{q^e-1}(a) = \begin{cases} 0 & \text{if $a\in\calT_0 \cup \{0\}$} \\
a^2 & \text{if $a\in\calT_1$,} \end{cases}
\qquad
D_{q^e+1}(a) = \begin{cases} a^2 & \text{if $a \in \calT_{0} \cup \{0\}$}\\
0 & \text{if $a\in\calT_1$.} \end{cases}
$$
\end{proposition}

\begin{proof}  For arbitrary $k\ge 1$ and any odd $d > 1$,
$D_k(\del d)=\del 1=\{0\}$ iff $d$ divides $k$
by Lemma~\ref{Dkdelta}.   Since $D_k(0)=0$ always and $D_0(x)=2=0$, 
it follows that for any $k\ge0$ and any odd $d>1$,
$D_k(\delstar d)=\{0\}$ iff $d$ divides $k$.  
By Proposition~\ref{calT_prop},
$\delstar {q-1}=\calT_{0}$ and $\delstar{q+1}=\calT_{1}$.
({\it i}) now follows.
For ({\it ii}), let $k=\pm 2 \pmod{q+1}$. Let
$b=\ang u\in \calT_{1}$, so $u^{q+1}=1$. Then
$D_k(b) = \ang{u^k}=\ang{u^{\pm2}}=\ang{u^2}=\ang u^2=b^2$.
Conversely, suppose $D_k(b)=b^2$ for all $b \in \calT_{1}$.
Let $u$ have order exactly $q+1$ and $b=\ang u$. Since $b\in\calT_{1}$,
$D_k(b)=b^2$.  Then $\ang{u^k}=\ang{u^2}$, so 
$u^k \in \{u^2,u^{-2}\}$ by~(\ref{angEquality}). Thus, $u^{k-2}=1$ or 
$u^{k+2}=1$.
It follows that $q+1$ divides $k-2$ or $k+2$, {\it i.e.} $k=\pm2 \pmod{q+1}$.
The second statement in ({\it ii}) is proved analogously.  Statement
({\it iii}) follows from ({\it i}) and ({\it ii})
by noting that
$$q^e= ((q-1)+1)^e \equiv 1 \pmod{q-1},\qquad
q^e= ((q+1)-1)^e \equiv (-1)^e \pmod{q+1}$$
and $D_k(a)=0=a^2$ when $a=0$.
\end{proof}

\section{Alternative descriptions of $\delstar{q-1}$ and $\delstar{q+1}$ in odd characteristic} \label{sec:oddq}

Proposition~\ref{calT_prop} provides an alternative description of the sets 
$\del{q-1}$ and $\del{q+1}$ in characteristic~2. 
We have found alternative descriptions in odd characteristic as well.
If $q$ is odd, define the following subsets of $\F_q$:
\begin{equation}\calZ = \{2,-2\},\qquad \calS_q = \{ a \in \F_q : \text{$a^2-4 $ is a nonzero square in $\F_q$}\},\label{ZSN} \end{equation}
$$\calN_q = \{ a \in \F_q : \text{$a^2-4$ is a nonsquare in $\F_q$}\}.$$
Note that $\F_q$ is a disjoint union:
\begin{equation} \F_q = \calS_q \sqcup \calN_q \sqcup \calZ.   \label{FqDecomp} \end{equation}
Proposition~\ref{prop:calT_prop2} below was observed by Brewer \cite{Brewer}. Theorem~\ref{thm:oddq} seems to be new.

\begin{proposition} \label{prop:calT_prop2} ($q$ odd.) 
$\calS_q = \delstar{q-1}$, $\calN_q = \delstar{q+1}$, and $|\calS_q| = (q-3)/2$, $|\calN_q|=(q-1)/2$.
\end{proposition}
\begin{proof}  Let $a \in\F_q$, and write $a = u + 1/u$, or equivalently $u^2 - a u + 1 = 0$. 
By the quadratic formula, $u=(a \pm \sqrt{a^2-4})/2$.  Thus, $u\in\F_q^\x$ iff $a^2-4$ is a square iff $a\in\calS_q\sqcup\calZ$.
In other words, $a=\ang u\in\calS_q\sqcup\calZ$ iff $u\in\mu_{q-1}$.
This shows that $\calS_q\sqcup\calZ=\del{q-1}=\delstar{q-1}\sqcup\calZ$, {\it i.e.}, $\calS_q=\delstar{q-1}$.
Now (\ref{FqDisjoint}) and (\ref{FqDecomp}) together imply $\calN_q=\delstar{q+1}$.
The last assertion follows from (\ref{deltaCount}). 
\end{proof}

\begin{theorem}  \label{thm:oddq} ($q$ odd.) 
For any $k\ge 1$, $D_k(\calZ) \subset \calZ$, $D_k(\calS_q) \subset \calS_q \sqcup \calZ$, and 
$D_k(\calN_q) \subset \calN_q \sqcup \calZ$. Further, $D_k$ permutes $\calZ$ iff $k$ is odd, and
it permutes $\calN_q$ iff $\GCD(k,q+1)=1$. If $\calS_q$ is nonempty (equivalently, $q>3$), then 
$D_k$ permutes $\calS_q$ iff $\GCD(k,q-1)=1$. If $(q-1)/2$ divides $k$ then $D_k(\F_q)\subset \calN_q \sqcup \calZ$.
If $(q+1)/2$ divides $k$ then $D_k(\F_q) \subset \calS_q \sqcup \calZ$.
\end{theorem}
\begin{proof} Since $2=\ang{1}$ and $-2=\ang{-1}$, 
\begin{equation} \text{$D_k(2)=\ang{1^k}=2$ and $D_k(-2)=\ang{(-1)^k} = (-1)^k 2$.} \label{DkZ} \end{equation}
This shows that $D_k(\calZ)
\subset \calZ$, and that $D_k$ permutes $\calZ$ exactly when $k$ is odd.
Next, since $D_k(\del d) \subset \del d$ for any~$d$, Proposition~\ref{prop:calT_prop2} implies that
$$D_k(\calS_q) \subset D_k(\calS_q \sqcup \calZ) = D_k(\del {q-1}) \subset \del {q-1} = \calS_q \sqcup \calZ$$
and
$$D_k(\calN_q) \subset D_k(\calN_q \sqcup \calZ) = D_k(\del {q+1}) \subset \del{q+1} = \calN_q \sqcup \calZ.$$
Since $\calN_q=\delstar{q+1}$, $D_k(x)$ permutes $\calN_q$ iff GCD$(k,q+1)=1$
by Corollary~\ref{cor:permuteDel}. Likewise,
since $\calS_q = \delstar{q-1}$, if $q-1>2$ then $D_k$ permutes
$\calS_q$ iff GCD$(k,q-1)=1$.
If $(q-1)/2$ divides~$k$, then by
Propositions~\ref{prop:dicksonImage} and~\ref{prop:calT_prop2},
$$D_k(\F_q) \subset \del2 \cup \del{q+1}
=\calZ \sqcup \delstar{q+1} = \calZ \sqcup \calN_q.$$
Similarly, if $(q+1)/2$ divides~$k$, then 
$D_k(\F_q) \subset \del{q-1} \cup \del{2}
=\calS_q \sqcup \calZ$.
\end{proof}

\begin{lemma} \label{lem:oddq}  ($q$ odd.) Let $k,\ell \ge 0$. If $\calS_q$ is nonempty, then
$D_k(a)=D_\ell(a)$ for all $a \in \calS_q$ iff $k\equiv \pm \ell \pmod{q-1}$. 
Also,  
$D_k(b)=D_\ell(b)$ for all $b \in \calN_q$ iff $k\equiv \pm \ell \pmod{q+1}$.
\end{lemma}
\begin{proof} In general, $D_k(a)=D_\ell(a)$ for all $a \in \del d$
iff $\ang{u^k}=\ang{u^\ell}$ for all $u\in \mu_d$ iff 
$u^k \in \{u^\ell,u^{-\ell}\}$ for all $u \in \mu_d$ iff 
$k \equiv \pm \ell\pmod d$. 
If $d>2$, then by taking $u$ primitive of order $d$, the above argument 
tells us that
$D_k(a)=D_\ell(a)$ for all $a \in \delstar d$ iff $k\equiv \pm\ell
\pmod d$. Now the result follows by noting that $\calS_q = \delstar{q-1}$
and $\calN_q = \delstar{q+1}$.
\end{proof}

\begin{corollary} ($q$ odd.) Let $e\ge0$. If $e$ is even, then $D_{q^e-1}(a)=2$ and $D_{q^e+1}(a)=
a^2-2$ for all $a \in \F_q$.  If $e$ is odd, then for $a \in \F_q$,
\begin{equation*}
D_{q^e-1}(a) = \begin{cases} a^2 - 2 & \text{if $a\in\calN_q$}\\
2 & \text{otherwise;} \end{cases}
\qquad
D_{q^e+1}(a) = \begin{cases} 2 & \text{if $a\in\calN_q$}\\
a^2-2 & \text{otherwise.} \end{cases}
\end{equation*}
\end{corollary}
\begin{proof}
Apply Lemma~\ref{lem:oddq} and eq.~(\ref{DkZ}),  noting that 
$q^e \equiv 1 \pmod {q-1}$,
$q^e \equiv (-1)^e \pmod{q+1}$, and $D_2(x)=x^2-2$.
Note that when $k$ is even, $D_k(a)=2=a^2-2$ for $a\in\{2,-2\}$.
\end{proof}

\section{Further refinements when $q$ is odd} \label{sec:refinements}
Theorem~\ref{thm:oddq} shows that $D_k$ preserves certain subsets of $\F_q$: it maps $\calZ$ to itself, $\calS_q$ to $\calS_q \sqcup \calZ$,
and $\calN_q$ to $\calN_q \sqcup \calZ$.  Further, precise conditions on $k$ are given that determine when $D_k$ permutes $\calZ$, $\calS_q$, and
$\calN_q$.
It turns out we can obtain finer results of this type when $q$ is odd, as we now describe.
Let $\jacobi aq$ denote the Legendre symbol, defined by (\ref{jacobiDef}).
It is well known that $\F_q^\x$ is cyclic; if $g$ is a generator then $g^i$ is a square in $\F_q^\x$ iff $i$ is even.
From this, we see that $\jacobi {g^k}q = (-1)^k$, and so the Legendre symbol is multiplicative:
$\jacobi a q \jacobi b q = \jacobi {ab}q$.
Define $\calA_\lambda^{\varepsilon_1,\varepsilon_2}$ by (\ref{Aij}). 
Noting that $\jacobi {a^2-4}q = \jacobi {a-2}q \jacobi {a+2}q$, we see that $\calS_q$ and $\calN_q$ break into disjoint unions:
\begin{eqnarray}
\calS_q &=& \calA_2^{++} \sqcup \calA_2^{--} = \{ a \in \F_q : \text{$a^2-4 $ is a nonzero square in $\F_q$}\},\label{ZSN2} \\
\calN_q &=& \calA_2^{+-} \sqcup \calA_2^{-+} = \{ a \in \F_q : \text{$a^2-4$ is a nonsquare in $\F_q$}\}. 
\end{eqnarray}
In addition,  $\F_q$ is a disjoint union:
\begin{equation} \F_q = \calA_2^{++} \sqcup \calA_2^{+-} \sqcup \calA_2^{-+} \sqcup \calA_2^{--} \sqcup \calZ = \calS_q \sqcup \calN_q \sqcup \calZ.   \label{FqDecomp2} \end{equation}

\begin{theorem} \label{thm:AijDescription} ($q$ odd.)
For each $\varepsilon_1,\varepsilon_2 \in \{1,-1\}$,
\begin{equation}
\calA_2^{\varepsilon_1,\varepsilon_2} = \left\{ v^2 + v^{-2} : v^{q-\varepsilon_1\varepsilon_2} = \varepsilon_2, v^4 \ne 1\right \}.
\label{AijCharacterization}
\end{equation}
Further, 
$$\calA_2^{++} = \delstar{(q-1)/2},\quad \calA_2^{--} = \delstar{q-1}\setminus \delstar{(q-1)/2},$$ 
$$\calA_2^{-+}=\delstar{(q+1)/2},\quad \calA_2^{+-}=\delstar{q+1}\setminus \delstar{(q+1)/2},$$
$$|\calA_2^{++}| = \lfloor (q-3)/4 \rfloor,\quad |\calA_2^{+-}| = \lfloor (q+1)/4 \rfloor,\quad 
|\calA_2^{-+}| = \lfloor (q-1)/4 \rfloor,\quad |\calA_2^{--}|=\lfloor (q-1)/4 \rfloor. $$
\end{theorem}
\begin{proof}   If $a=\ang{v^2}$, where $v\in\cj\F_q$, then by Lemma~\ref{lem:FqDecomposition},
$$a \in \F_q \iff v^2 \in \mu_{q-1} \cup \mu_{q+1} \iff v \in \mu_{2q-2}\cup\mu_{2q+2}.$$
Now $v^{2q-2}=1 \iff v^{q-1}\in\{1,-1\}$, and similarly $v^{2q+2}=1 \iff v^{q+1}\in\{1,-1\}$.  
Thus, there are integers $c,d\in \{1,-1\}$ such that $v^{q-c}=d$.
Note that $v^4=1$ iff $\ang{v^2} \in \ang{\mu_2}=\calZ$. Assume now that $v^4\ne 1$, therefore $\ang{v^2}\not \in \calZ$. 
Then the integers 
$c,d\in \{1,-1\}$ such that $v^{q-c}=d$ are uniquely determined. Define 
$\typ(v)=(c,d)$. Let $a=\ang{v^2}\in\F_q$. Then 
$$a - 2 = (v-1/v)^2, \qquad a + 2 = (v+1/v)^2 $$
and these are nonzero because of the assumption that $v^4 \ne 1$. Now $a-2$ is a square iff $v-1/v \in \F_q$.
In general, if $v,w$ are nonzero then $v-1/v=w-1/w$ iff $w \in \{v,-1/v\}$. Thus, $(v-1/v)=(v-1/v)^q$ iff $v^q \in \{v,-1/v\}$,
which is equivalent to $\typ(v) \in \{(1,1),(-1,-1)\}$.  It follows that if $\typ(v)=(c,d)$ then $\jacobi{a-2}q = cd$.
Similarly, $a+2$ is a square iff $v+1/v \in \F_q$
iff $v^q \in \{v,1/v\}$ iff $\typ(v) \in \{(1,1),(-1,1)\}$. It follows that $\jacobi{a+2}q=d$.
This proves (\ref{AijCharacterization}). In particular,
\begin{eqnarray*}
\calA_2^{++} &=& \{ \ang{v^2} : v^{q-1}=1,\ v^4\ne 1 \} \\
&=&  \{\ang u : u^{(q-1)/2} = 1,\ u^2 \ne 1 \} = \del{(q-1)/2}\setminus \del 2  \\
&=& \delstar {(q-1)/2},
\end{eqnarray*}
and $\calA_2^{-+} = \delstar {(q+1)/2}$ by a similar computation.  Next,
 $\calA_2^{--} = \calS_q \setminus \calA_2^{++} = \delstar{q-1} \setminus \delstar{(q-1)/2}$, and $\calA_2^{+-}=\calN_q\setminus \calA_2^{-+}
=\delstar{q+1} \setminus \delstar{(q+1)/2}$. 
By (\ref{deltaCard}), $|\calA_2^{++}| =  \lfloor (d-1)/2 \rfloor$, 
where $d=(q-1)/2$; this proves $|\calA_2^{++}|= \lfloor (q-3)/4 \rfloor$. The same computation with $d=(q+1)/2$ shows $|\calA_2^{-+}|=
\lfloor (q-1)/4 \rfloor$.  Next, $|\calA_2^{--}| = |\delstar{q-1}|-|\delstar{(q-1)/2}| = (q-3)/2 - \lfloor (q-3)/4 \rfloor = \lceil (q-3)/4 \rceil$.
In general, if $2m$ is an integer, then $\lceil m\rceil=\lfloor m+1/2 \rfloor$.
Thus, $\lceil{(q-3)/4}\rceil = \lfloor  (q-1)/4 \rfloor$. The formula $|\calA_2^{+-}|=\lfloor (q+1)/4 \rfloor$ may be computed similarly.
\end{proof}

\begin{lemma} \label{lemma:Z01}  ($q$ odd.)
$\calZ \cap \del{(q-1)/2} = \{2,-2\varepsilon\}$
and $\calZ \cap \del{(q+1)/2} = \{2,2\varepsilon\}$, 
where $\varepsilon=(-1)^{(q-1)/2}$. 
\end{lemma}
\begin{proof} In general, $\mu_d\cap\calZ=\mu_d\cap\mu_2=\mu_g$, where
$g=\GCD(d,2)$. Thus, $\mu_d\cap\calZ=\calZ$ if $d$ is even, 
and $\mu_d\cap\calZ=\mu_1=\{2\}$
if $d$ is odd. Noting that $(q-\varepsilon)/2$ is even and $(q+\varepsilon)/2$
is odd, the result follows.
\end{proof}

Because of Theorem~\ref{thm:AijDescription}, $D_k$ maps the sets $\calA_2^{\varepsilon_1,\varepsilon_2}$ to one another in patterned ways, as described in 
Theorem~\ref{thm:Aij}.  We are now ready  to prove that theorem.

\noindent {\it Proof of Theorem~\ref{thm:Aij}.}\qquad 
If $\typ(v)=(c,d)$, {\it i.e.}, $v^{q-c}=d$ and $v^4\ne 1$, then $(v^k)^{q-c}=d^k$.   
Thus, either $v^{4k} = 1$ or $\typ(v^k)=(c,d^k)$.
By (\ref{AijCharacterization}),
$\typ(v)=(c,d)$ iff $\ang{v^2} \in \calA_2^{cd,d}$,  therefore
$$D_k(\calA_2^{cd,d}) \subset \calA_2^{cd',d'} \cup \{2,-2\},\quad\text{where $d'=d^k$.}$$
To complete the proof of ({\it i})--({\it iii}), we need only consider 
$D_k(\calA_2^{\varepsilon_1,\varepsilon_2})\cap \calZ$.
First, $D_k(\calA_2^{++}) \subset D_k( \del{(q-1)/2}) \subset \del{(q-1)/2}$,
so $D_k(\calA_2^{++})\cap\calZ \subset \{2,-2\varepsilon\}$ by 
Lemma~\ref{lemma:Z01}.
Similarly, $D_k(\calA_2^{-+}) \subset \del{(q+1)/2}$, so 
$D_k(\calA_2^{-+})\cap\calZ \subset \{2,2\varepsilon\}$ by Lemma~\ref{lemma:Z01}. 
This proves ({\it i}).
For ({\it ii}), let $\ang{v^2} \in \calA_2^{-c,-}$, so $\typ(v)=(c,-1)$, and let $k$ be odd. Then $(v^k)^{q-c}=-1$. 
If $c=\varepsilon$, then
$(v^k)^{q-c}$ is a power of $v^{4k}$, therefore $v^{4k}\ne 1$, showing $D_k(\ang{v^2}) \not \in\calZ$.  
Thus, $D_k(\calA_2^{-\varepsilon,-}) \cap \calZ = \emptyset$. Next, if 
$c=-\varepsilon$ then
$(v^k)^{q-c}$ is an odd power of $v^{2k}$. Since $(v^k)^{q-c}=-1$, we see that $v^{2k}\ne 1$, and therefore $D_k(\ang{v^2})\ne 2$. Thus,
$D_k(\calA_2^{-\varepsilon,-}) \cap \calZ \subset \{-2\}$. 
To prove ({\it iii}), let $k$ be even. 
Then 
$D_k(\calA_2^{--}) \subset D_k(\del{q-1}) \subset \del{(q-1)/2}$,
so $D_k(\calA_2^{--})\cap\calZ \subset \{2,-2\varepsilon\}$.  Likewise,
$D_k(\calA_2^{+-}) \subset D_k(\del{q+1}) \subset \del{(q+1)/2}$, so
$D_k(\calA_2^{+-})\cap\calZ \subset \{2,2\varepsilon\}$.

For ({\it iv}), if $(q-1)/2$ divides $k$ then $D_k(\del{q-1}) \subset \del 2 = \calZ$, so
$D_k(\F_q) \subset D_k(\del{q+1}) \cup \calZ \subset \del{q+1}  \cup \calZ = \calA_2^{+-} \cup \calA_2^{-+} \cup \calZ$.
If in addition $k$ is even, then $D_k(\F_q)\subset D_k(\del{q+1})\cup\calZ \subset \del{(q+1)/2}\cup\calZ=
\calA_2^{-+} \cup\calZ$.
The statement ({\it v}) is proved similarly.

Finally, we prove ({\it vi}). In general, if $p\nmid d$ and $\GCD(k,d)=1$ then $D_k$ permutes $\del d$. Since $D_k$ takes $\calZ$ to itself,
$D_k$ must also permute $\delstar d$.  If in addition $d'$ divides $d$, then since $D_k$ takes $\delstar{d'}$ to itself, we see that 
$D_k$ must permute $\delstar{d} \setminus \delstar{d'}$. These considerations show that $D_k$ permutes $\calA_2^{\varepsilon_1,\varepsilon_2}$
when $\calA_2^{\varep_1,\varep_2}$ is nonempty and GCD$(k,d^{\varepsilon_1,\varepsilon_2})=1$. Now we prove the converse. Since $\calA_2^{++}=\delstar{(q-1)/2}$, 
if this set is nonempty then there is $\ang b \in \calA_2^{++}$ with order$(b)=(q-1)/2$.
If $\GCD(k,(q-1)/2)) > 1$, then every element $\ang {u} \in D_k(\calA_2^{++})$ has order$(u) < (q-1)/2$,
showing that $\ang b$ is not in the image.  Thus, $D_k$ permutes $\calA_2^{++}$
if and only if $\GCD(k,(q-1)/2)=1$. The other cases are proved similarly.
\qed

\medskip
\noindent{\bf Example.} 
In $\F_{29}$, the nonzero squares are $\{1,4,5,6,7,9,13,16,20,22,23,24,25,28\}$.
By the definition of $\calA_2^{\varepsilon_1,\varepsilon_2}$, 
$$\calA_2^{++}=\{3,7,11,18,22,26\}, \qquad \calA_2^{+-}=\{1,6,8,9,15,24,25\},$$
$$\calA_2^{-+}=\{4,5,14,20,21,23,28\},  \qquad \calA_2^{--} = \{\, 0, 10, 12, 13, 16, 17, 19 \,\}. $$
We will demonstrate Theorem~\ref{thm:Aij} in some special cases.
One assertion is that if $k$ is odd, then $D_k(\calA_2^{--})\subset \calA_2^{--}$, and 
$D_k$ permutes $\calA_2^{--}$ if
and only if $(k,29-1)=1$. We test this in a few cases. 
For $k=5$, $D_5(x) = x^5-5 x^3 + 5x$. 
$D_5(x)$ permutes $\calA_2^{--}$ as follows: it takes $0,10,12,13,16,17,19$ to $0,17,16,19,10,13,12$, respectively. 
This is the permutation $(0)(10,17,13,19,12,16)$.
On the other hand, $D_7=x^7 - 7 x^5 + 14 x^3 - 7 x$ maps $\calA_2^{--}$ to $\{0\}$, so $D_7$ maps $\calA_2^{--}$ to itself but is not a permutation. 
In this case, the seven elements of $\calA_2^{--}$ are roots of $D_7$, so
$$D_7(x)=\prod_{a \in \calA_2^{--}} (x-a) \pmod{29}.$$ 

The above formula for $D_7$ is quite interesting, and the phenomenon generalizes: If
$q\equiv 1 \pmod 4$ then $D_{(q-1)/4}(x) = \prod_{a \in \calA_2^{--}} (x-a)$. This and 
similar results 
will be shown in the next section.

\section{Wilson-like theorems} \label{sec:wilson-like}

In this section, we use Dickson polynomials to obtain some Wilson-like theorems and other results in elementary
number theory.  The starting point is
Theorem~\ref{thm:fijDE} below, which generalizes the phenomenon that was observed for $D_7(x)$ at the end of the preceding section. 
First we need to introduce more sets and polynomials.

For $\varepsilon_1,\varepsilon_2 \in \{1,-1\}$ and $\lambda \in \F_q^\x$, 
define for $q$ odd:
\begin{equation} \calB_\lambda^{\varepsilon_1,\varepsilon_2} = 
\left\{ b \in \F_q : \jacobi{\lambda-b}q=\varepsilon_1,\
\jacobi{\lambda+b}q=\varepsilon_2 \right\}. \label{BDef}
\end{equation}
Note that $\calB_\lambda^{\varepsilon_1,\varepsilon_2} = \calA_\lambda^{\varepsilon\varepsilon_1,\varepsilon_2}$,
where as usual $\varepsilon = \jacobi{-1}q=(-1)^{(q-1)/2}$. Also note that $b \mapsto -b$ gives
a bijection from $\calB_\lambda^{\varepsilon_1,\varepsilon_2}$ to $\calB_\lambda^{\varepsilon_2,\varepsilon_1}$.
As usual, we often write the superscripts as $+$ and $-$ instead of 1 and $-1$.
Let $m = (q-\varep)/4$, which is always an integer. By Theorem~\ref{thm:AijDescription},
\begin{eqnarray}  
\calB_2^{++} &=& \calA_2^{\varep,+} = \delstar{2m}, \qquad
\calB_2^{--} = \calA_2^{-\varep,-} = \delstar{4m} \setminus \delstar{2m},  \label{BppBmm} \\
\calB_2^{-+} &=& \calA_2^{-\varep,+} = \delstar{2m+\varep}, \qquad
\calB_2^{+-} = \calA_2^{\varep,-} = \delstar{4m+2\varep} \setminus \delstar{2m+\varep}.  \label{BpmBmp} 
\end{eqnarray}

Define the following polynomials in $\F_q[x]$ for $\varep_1,\varep_2 \in \{1,-1\}$:
\begin{equation} \label{fijDef}
f^{\varepsilon_1,\varepsilon_2}(x) = \prod_{a \in \calA_2^{\varepsilon_1,\varepsilon_2}} (x-a), \qquad 
g^{\varepsilon_1,\varepsilon_2}(x) = \prod_{a \in \calB_2^{\varepsilon_1,\varepsilon_2}} (x-a) = f^{\varep\varep_1,\varep_2}(x)
\end{equation}
\begin{equation} \label{fSfNdef}
f_\calS = \prod_{a\in\calS_q} (x-a), \qquad
 f_\calN = \prod_{a\in\calN_q} (x-a), \qquad
 f_\calZ = \prod_{a \in \calZ}(x-a) = x^2-4. 
\end{equation}
Since $\calS_q = \calA_2^{++} \sqcup \calA_2^{--}$ and $\calN_q = \calA_2^{+-}\sqcup \calA_2^{-+}$,
$$ f_\calS = f^{++} f^{--},\qquad f_\calN = f^{+-} f^{-+}.$$
Also, since $\prod_{a\in\F_q}(x-a) = x^q-x$, (\ref{FqDecomp}) implies that 
$$f_\calS(x) f_\calN(x) f_\calZ(x) = x^q-x.$$

Define $E_k\in\Z[x]$ for $k\ge 0$ by 
\begin{equation} \label{Edef} \text{$E_0(x)=1$, $E_1(x)=x$, and 
$E_{k+2}(x) = x E_{k+1}(x) - E_k(x)$.} \end{equation}
This is called a Dickson polynomial of the second kind, and like $D_k$, it has
been widely studied. It is well known (and can easily be shown by 
induction on $k$) that 
\begin{equation} \label{Efunctional} E_{k-1}(\ang u) = \frac{u^k - 1/u^k}{u-1/u}
\end{equation}
Parts {(\it i}) and {(\it ii)} of the next theorem are known. (See \cite{chou1,BZ}).

\begin{theorem} \label{thm:fijDE} ($q$ odd.)  In $\F_q[x]$, the following holds. \\
({\it i})\quad If $k>1$, then 
\begin{equation}
E_{k-1}(x) = \prod_{a \in \delstar{2k}} (x-a). \label{Edelta}
\end{equation}
({\it ii})\quad If $k\ge 1$, then 
\begin{equation} 
D_k(x) = \prod_{ a \in \delstar {4k}\setminus \delstar {2k}} (x-a).
\label{Ddelta}
\end{equation}
({\it iii})\quad $f_\calS(x) = E_{(q-3)/2}(x)$ and $f_\calN(x) = E_{(q-1)/2}(x)$. \\
({\it iv})\quad Let $m=(q-\varep)/4$. Then $D_{m}(x) = g^{--}(x)$
and $E_{m-1}(x) = g^{++}(x)$, where $g^{\pm,\pm}(x)$ is defined in~(\ref{fijDef}). \\
\end{theorem}

\begin{proof}  First we prove the assertions about $E_k$.
The two polynomials in (\ref{Edelta}) have the same degree by (\ref{deltaCard}), 
so to prove equality it suffices
to show that if $a = \ang u \in \delstar{2k}$ then $E_{k-1}(a) = 0$.
Here $u \in \mu_{2k}$ and $u^2\ne 1$. Since $u^{2k}=1$, $u^k = 1/u^k$.
On the other hand, $u \ne 1/u$. Thus, $E_{k-1}(\ang u)=0$ by (\ref{Efunctional}), as
required.   All the assertions about $E_k$ follow from this observation, together
with the fact that $\calS_q = \delstar{q-1}$, $\calN_q = \delstar{q+1}$,
and $\calB_2^{++} = \delstar{(q-\varep)/2} = \delstar{2m}$. 

It remains to prove the assertions about $D_k$. 
The polynomials in (\ref{Ddelta}) are monic and have the same degree by
(\ref{deltaCard}), so to show equality, it suffices to show that $a \in \delstar{4k} \setminus
\delstar{2k}$ implies $D_k(a)=0$. Indeed, writing $a=\ang u$ we have
$u^{4k}=1$ but $u^{2k}\ne 1$, so $u^{2k}=-1$ and $u^k=i$ is a square root of
$-1$. Then $D_k(a)=\ang{u^k}=\ang{i}=0$, as required. This proves ({\it ii}), and ({\it iv}) follows
immediately by (\ref{BppBmm}).
\end{proof}

We are ready to prove the Wilson-like theorems.  For many more theorems of this sort, see the sequel to this article, \cite{wilson-like}.

\begin{theorem} \label{thm:wilson-like}
($q$ odd.) ({\it i})\ $\prod \{\,a\in\F_q^\x : \text{$a$ and $4-a$ are nonsquares} \,\}=2$. \\
({\it ii})\  $\prod \{\,a \in \F_q^\x : \text{$-a$ and $4+a$ are nonsquares} \,\} = \begin{cases} 2 & \text{if $q\equiv \pm1\pmod8$}
\\ - 2 & \text{if $q\equiv \pm3 \pmod8$.}\end{cases}$
\end{theorem}

\begin{proof} Let $m=(q-\varep)/4$. Substituting $x=2$ into the formula $D_m(x) = g^{--}(x)$ gives
$$\prod \{ (2-b) : b \in \calB_2^{--} \} = 2,$$
since $D_n(2)=2$ for all $n$. Set $a=2-b$. Then $2+b=4-a$. The condition that $b \in \calB_2^{--}$ is
therefore equivalent to $a$ and $4-a$ being nonsquares.  This proves~({\it i}).

For the second statement, evaluate $D_m(x) = g^{--}(x)$ at $x=-2$. Since $D_n(-2)=D_n(\ang{-1})=\ang{(-1)^n} = (-1)^n \cdot 2$ for all~$n$,
$$\prod \{ (-2-b) : b \in \calB_2^{--} \} = (-1)^{m}2.$$
Set $a=-2-b$. Then $2+b=-a$ and $2-b=4+a$, so $b\in\calB_2^{--}$ iff $-a$ and $4+a$ are nonsquares.
This proves~({\it ii}).
\end{proof}

\section{More number-theoretic results} \label{sec:numberTheory}

The statements of the next few propositions do not mention 
Dickson polynomials, although they are certainly working behind the scenes! Let $\calB_\lambda^{\varep_1,\varep_2}$ be as in (\ref{BDef}).

\begin{proposition} \label{prop:D2}  ($q$ odd.) The map $b \mapsto b^2-2$ gives a permutation $\pi$ on $\calB_2^{-+}$, and the inverse permutation is
$$ \pi^{-1}(b) = \prod\{b-a : a \in \calB_2^{--} \}\qquad\text{for $b \in\calB_2^{-+}$.}$$
Further,  $\pi^{-1}(b)$ is the unique square root of $2+b$ such that $2+\sqrt{2+b}$ is a square. \\
\end{proposition}

\begin{proof}  Noting that $\calB_2^{-+}=\delstar{2m+\varep}$ by (\ref{BpmBmp}), both $D_2$ and $D_m$ permute it by Corollary~\ref{cor:permuteDel}. 
Further, if $b=\ang u \in \calB_2^{-+}$ then $D_2(b)=\ang{u^2}=\ang{u}^2-2=b^2-2=\pi(b)$, and
$$D_m \circ \pi(b) = D_m(\ang{u^2})=\ang{u^{2m}}=\ang{u^{-\varep}}=\ang u=b.$$
Thus, $D_m=\pi^{-1}$.  By Theorem~\ref{thm:fijDE}({\it iv}), $D_m(b)=\prod\{b-a : a \in \calB_2^{--}\}$. 
Let $c=\pi^{-1}(b)$. Then $b=c^2-2$, so $c$ is a square root of $b+2$. Noting that $c\in\calB_2^{-+}$ and
$-c \in \calB_2^{+-}$, the choice of square root is determined by the property that $\sqrt{b+2}+2$ is a square.
\end{proof}

\begin{proposition} ($q$ odd.)  \label{prop:fSfNvalues} Let $c \in \F_q$. Then
	$$E_{(q-3)/2}(c)=\prod \left\{ c-a : a \in \F_q,\ \jacobi{a^2-4}q=1\right\} = \begin{cases} 
0 & \text{if $\jacobi{c^2-4}q=1$} \\
\jacobi{c-2} q & \text{if $\jacobi{c^2-4}q=-1$} \\
-1/2 & \text{if $c=2$ } \\
\jacobi{-1}q/2 & \text{if $c=-2$ }
\end{cases}
$$
	$$E_{(q-1)/2}(c) = \prod \left\{ c-b : b \in \F_q,\ \jacobi{b^2-4}q=-1\right\} = \begin{cases} 
\jacobi{c-2} q & \text{if $\jacobi{c^2-4}q=1$} \\
0 & \text{if $\jacobi{c^2-4}q=-1$} \\
1/2 & \text{if $c=2$ } \\
\jacobi{-1}q/2 & \text{if $c=-2$. }
\end{cases}
$$
\end{proposition}

\begin{proof} We prove the first formula only, as the second is similar. 
By Theorem~\ref{thm:fijDE}, the product equals $E_{(q-3)/2}(c)$.
	If $q=3$ the product is an empty product, which by convention equals~1.
The expression on the right is 1 for each $c\in \{0,1,2\}=\F_3$.  Thus, the two sides agree if $q=3$. Now suppose that $q>3$. 
If $\jacobi{c^2-4}q=1$ then the product vanishes since it contains a factor $(c-c)$. If $\jacobi{c^2-4}q=-1$ then $c \in \calA_2^{+-} \sqcup \calA_2^{-+}$
and consequently $c = \ang{w^2}$, where $w^{q+1}=\jacobi{c+2}q$ and $w^4\ne 1$, by Theorem~\ref{thm:AijDescription}. Then
\begin{eqnarray*}
E_{(q-3)/2}(c) &=& (w^{q-1}-w^{1-q})/(w^2-w^{-2}) \\
&=&  ( w^{q+1-2} - w^{-(q+1-2)} )/(w^2-w^{-2}) \\
&=& \jacobi{c+2}q (w^{-2}-w^2)/(w^2-w^{-2})=-\jacobi{c+2}q=\jacobi{c-2}q.
\end{eqnarray*}

For any $k$,
$$E_{k-1}(\ang u)=\frac{u^k-u^{-k}}{u-u^{-1}} = u^{1-k} \frac{u^{2k}-1}{u^2-1}.$$
By expanding the right side as $u^{1-k}(1+u^2+u^4 + \cdots u^{2k-2})$ and then evaluating at $u = \pm1$, we obtain that $E_{k-1}(2)=k$
and $E_{k-1}(-2)=(-1)^{k-1}k$.  Applying this with $k-1=(q-3)/2$, and observing that $(-1)^k=\jacobi{-1}q$, we obtain the formula when $c=\pm2$.
\end{proof}

\begin{corollary} ($q$ odd.) \label{cor:fSfNvalues} For any $c \in \F_q$, 
the equation (\ref{jacc}) holds.
\end{corollary}
\begin{proof} In Proposition~\ref{prop:fSfNvalues}, the sum of the two right-hand sides
is $\jacobi{c-2}q$ if $c^2-4\ne 0$, $0$ if $c=2$,
and $\jacobi{-1}q$ if $c=-2$. In each case, this sum
is $\jacobi{c-2}q$. If we shift $c$ by~2, and also shift $a$ and $b$ by 2 in the left-hand sides,
then the formula~(\ref{jacc}) follows.
\end{proof}

\begin{corollary} \label{cor:fSfNvalues2} If $q=p^n$ is odd,
then $E_{(q-3)/2}(x) + E_{(q-1)/2}(x) \equiv (x-2)^{(q-1)/2} \pmod p$.
\end{corollary}
\begin{proof}  In Proposition~\ref{prop:fSfNvalues}, the sum of the two left sides is
$E_{(q-3)/2}(c)+E_{(q-1)/2}(c)$, and the sum of the two
right sides  is $(c-2)^{(q-1)/2}$. Thus,
$E_{(q-3)/2}(x)+E_{(q-1)/2}(x)-(x-2)^{(q-1)/2}$ has $q$ roots in $\F_q$.
Since its degree is less than $q$, the polynomial must vanish identically in $\F_q[x]$.
Since the polynomial has integer coefficients, this is equivalent to vanishing mod~$p$.
\end{proof}

Our next goal is to find 
analogues to the above results when $q$ is even. 
Let $q=2^n$. Part~({\it i}) of the following proposition is known (see \cite{chou1,BZ}).

\begin{proposition}   \label{prop:Dqpm1Factor} ($q$ even.) Let $\calT_{i}$ be as in~(\ref{TqjDef}) for $i\in\{0,1\}$, and consider
$D_k(x)$ as a polynomial in $\F_q[x]$. \\
{\it (i)}\  If $k$ is odd, then $D_k(x) = x f(x)^2$, where $f(x) = \prod_{a \in \delstar k} (x-a)$. \\
{\it (ii)}\ $D_{q+1}(x) = x f_1(x)^2$ and $D_{q-1}(x) = x f_0(x)^2$, where 
$f_i(x) = \prod_{a \in \calT_{i} } (x-a)$.
\end{proposition}
\begin{proof} It is well known (or can be easily shown by induction) that if $k$ is odd, then the coefficients
of $D_k(x)$ are zero for even powers of $x$. Thus, $D_k(x) = x f(x)^2$, where $f \in \F_2[x]$ and $f$ is monic.
Since $D_k(\delstar k)=\{0\}$, we see that every element of $\delstar k$ is a root of $f$.
Since $\deg(f) = (k-1)/2 = |\delstar k|$, we conclude that $f(x) = \prod_{a\in \delstar k} (x-a)$.
This proves~({\it i}), and ({\it ii}) now follows from~(\ref{Tq01}).
\end{proof}

\begin{proposition} \label{prop:sqrtc} ($q$ even.) Let $\calT_{0}$ and $\calT_{1}$ be as in (\ref{TqjDef}).
For all $c \in \F_q$,
$$\prod_{a \in \calT_{0}} (c+a) + \prod_{b \in \calT_{1}} (c+b) = c^{1/2}.$$ 
\end{proposition}
\begin{proof}  
By Proposition~\ref{prop:Dqpm1Factor}, $D_{q+1}(x) = x \prod_b (x+b)^2$ and
$D_{q-1}(x) = x \prod_a (x+a)^2$, where $b$ runs through $\calT_1$ and $a$ runs through $\calT_0$.  
On the other hand, if $c_0\in \calT_{0}$ 
and $c_1\in \calT_{1}$ then $D_{q+1}(c_0)=c_0^2$ and $D_{q-1}(c_1)=c_1^2$ by Proposition~\ref{prop:evenPermute}.
Thus, $c_0 \prod_b(c_0+b)^2 = c_0^2$ and $c_1 \prod_a(c_1+a)^2=c_1^2$.
It follows that $\prod_b (c_0+b)^2=c_0$, {\it i.e.} $\prod_b (c_0+b) = \sqrt{c_0}$,
and similarly $\prod_a (c_1+a) = \sqrt {c_1}$. 

Let $c\in\F_q$. 
Since $\F_q = \calT_{0} \sqcup \calT_{1} \sqcup \{0\}$, it suffices to prove the result in the cases
$c \in \calT_{0}$, $c\in \calT_{1}$, and $c=0$.

If $c \in \calT_{0}$, then $\prod_a (c+a) = 0 $ and $\prod_b (c+b) = c^{1/2}$. 
If $c \in \calT_{1}$, then $\prod_a (c+a) = c^{1/2} $ and $\prod_b (c+b) = 0$.
Thus, the proposition holds in those two cases.

Finally, if $c=0$ then we claim that $\prod_a (0+a) = \prod_b (0+b)=1$, so that again the proposition holds.
We have $\prod_a (1/a) = \prod \{v\in \F_q^\x : \Tr(v) = 0 \}$. Every such $v$ can be written
in exactly two ways as $r^2+r$, where $r \not \in \F_2$, so
$\prod_a (1/a)^2 = \prod_{r \in \F_q \setminus \F_2 } (r^2+r) = \prod_{r \in\F_q\setminus \F_2} r(r+1)$.
We break this as $\prod r \x \prod (r+1) = \left(\prod r\right)^2$, where these products are over $\F_q \setminus \F_2$.
By Wilson's Theorem, this equals $1$. We showed $\prod_a (1/a)^2=1$, therefore $\prod a = 1$. Since
$\calT_{0}\sqcup \calT_{1} = \F_q^\x$ and $\prod a \prod b = \prod F_q^\x = 1$,
$\prod b = 1$ as well.
\end{proof}

\section{Elementary symmetric functions of some natural sets in $\F_q$.}
\label{sec:sigmak}

Wilson-like formulas such as (\ref{T40}), (\ref{T04}) give closed formulas for the products of certain natural subsets of $\F_q^\x$, 
but one might ask about other symmetric functions.  The next theorem answers that question for four specific subsets.
If $B\subset\F_q$, let $\sigma_k(B)\in\F_q$ denote the elementary symmetric function of degree~$k$ in $B$, determined by the formula in $\F_q[x]$:
$$\prod_{b\in B} (x+b)= \sum_{k=0}^{|B|} \s_k(B) x^{|B|-k}.$$ 

\begin{theorem}  \label{thm:sigmak}
	Let $q$ be odd and $m=\lfloor(q+1)/4\rfloor$. \\
{\it (i)}\ Let $\calB_2^{--}=\set{a \in \F_q^\x : \jacobi{2-b}q=\jacobi{2+b}q=-1}.$ Then $|\calB_2^{--}|=m$, 
and for $0<k\le m$,
$$\s_k(\calB_2^{--}) = \begin{cases} 0 & \text{if $k$ is odd} \\ (-1)^i \left\{ \choose{m-i}i + \choose{m-i-1}{i-1} \right\} & \text{if $k=2i$ is even.} \end{cases}$$
{\it (ii)}\ Let $\calB_2^{++}=\set{a \in \F_q^\x : \jacobi{2-b}q=\jacobi{2+b}q=1}$. Then $|\calB_2^{++}|=m-1$, 
and for $0<k\le m-1$,
$$\s_k(\calB_2^{++}) = \begin{cases} 0 & \text{if $k$ is odd} \\ (-1)^i \choose{m-1-i}{i} & \text{if $k=2i$ is even.} \end{cases}$$
{\it (iii)}\ Let $\calS_q=\set{a \in \F_q^\x : \jacobi{a^2-4}q=1}$. Then $|\calS_q|=(q-3)/2$, and for $0<k\le (q-3)/2$,
$$\s_k(\calS_q) = \begin{cases} 0 & \text{if $k$ is odd} \\ (-1)^i \choose{(q-3)/2-i}{i} & \text{if $k=2i$ is even.} \end{cases}$$
{\it (iv)}\ Let $\calN_q=\set{a \in \F_q^\x : \jacobi{a^2-4}q=-1}$. Then $|\calN_q|=(q-1)/2$, and for $0<k\le (q-1)/2$,
$$\s_k(\calN_q) = \begin{cases} 0 & \text{if $k$ is odd} \\ (-1)^i \choose{(q-1)/2-i}{i} & \text{if $k=2i$ is even.} \end{cases}$$
	{\it (v)} For each $0<k\le m-1$, 
	$$\s_k(\calB_2^{++}) = 
	(1+(-1)^{(q-1)/2}4k) \s_k(\calB_2^{--}).$$
	{\it (vi)} For each $0<k\le (q-3)/2$, 
	$$(k+1)\s_k(\calS_q)=(2k+1)\s_k(\calN_q).$$
\end{theorem}

\begin{proof}
	{\it (i)}.\ By Theorem~\ref{thm:fijDE}({\it iv}),  $|\calB_2^{--}|=m$ and $\prod_{b\in \calB_2^{--}} (x-b) = D_m(x)$. 
	Since $\calB_2^{--}$ is permuted under $b\mapsto -b$, $\prod_{b\in\calB_2^{--}}(x-b)=\prod_{b\in\calB_2^{--}}(x+b)$. Thus,
	\begin{equation} \prod_{b\in \calB_2^{--}} (x+b) 
	=\sum_{0\le k\le m}\s_k(\calB_2^{--}) x^{m-k} =  D_m(x). \label{eq:B2mm} \end{equation}
	By \cite[Definition 1.1]{Lidl}, together with the identity 
	$$\frac n{n-i} \choose{n-i}i = \choose{n-i}i + \choose{n-1-i}{i-1}\qquad\text{for $0<i\le \lfloor n/2 \rfloor$},$$
	we have for $n>0$
	$$D_n(x) = x^n + \sum_{1\le i \le \lfloor n/2\rfloor } \left\{ \left( {n-i} \atop i \right)+\left({n-i-1} \atop {i-1}\right)\right\} (-1)^i x^{n-2i}.$$
	This, together with (\ref{eq:B2mm}), implies {\it (i)}.

	{\it (ii)--(iv)}.\	Each set
	$\calB_2^{++}$, $\calS_q$, $\calN_q$ is permuted under $b\mapsto -b$, therefore $\prod_b (x-b)=\prod_b(x+b)$ for these sets.
	Then by Theorem~\ref{thm:fijDE}{\it (iii)}, {\it (iv)}, 
	\begin{equation}  E_{m-1}(x) = \prod_{b\in\calB_2^{++}} (x+b)=\sum_{k=0}^{m-1} \s_k(\calB_2^{++}) x^{m-1-k}, \label{eq:B2pp} \end{equation}
$$E_{(q-3)/2}(x) = \prod_{b\in\calS_q} (x+b) 
	= \sum_{k=0}^{(q-3)/2} \s_k(\calS_q) x^{(q-3)/2-k}, $$
$$E_{(q-1)/2}(x) = \prod_{b\in\calN_q} (x+b) 
	= \sum_{k=0}^{(q-1)/2} \s_k(\calN_q) x^{(q-1)/2-k}.  $$
By \cite[Definition~2.2]{Lidl}, 
$$E_n(x) = \sum_{i=0}^{\lfloor n/2 \rfloor} \left( {n-i} \atop i \right) (-1)^i x^{n-2i}\in\Z[x].$$
This gives the formulas for $\sigma_k(\calB_2^{++})$, 
	$\sigma_k(\calS_q)$, and $\sigma_k(\calN_q)$ as claimed.

{\it (v)}. Schur \cite{Schur} showed that $E_{n-1}(x)=(1/n)D_n'(x)$
	for all $n\ge1$. (This will be explained in the next section; see
	(\ref{eq:DicksonDeriv}).) 
	Noting that $m=(q-\varep)/4\equiv -\varep/4 \pmod p$,
	where $\varep = (-1)^{(q-1)/2}$, the result is immediate
	from (\ref{eq:B2mm}) and (\ref{eq:B2pp}).

	{\it (vi)}. For any $0<i\le \lfloor n/2 \rfloor$, 
	$$\choose{n+1-i}i = \frac{n+1-i}{n+1-2i} \choose{n-i}i.$$
	Setting $n=(q-3)/2$ gives
	$$((q-1)/2-2i)\choose{(q-1)/2-i}i = ( (q-1)/2-i)\choose{(q-3)/2-i}i$$
	and so in $\F_q$,
	$$(-1-4i)\s_{2i}(\calN_q)=(-1-2i)\s_{2i}(\calS_q).$$
	This proves ({\it vi}) when $k$ is even. When $k$ is odd, the equality is trivial since both sides are zero.
\end{proof}

\section{Chebyshev polynomials} \label{sec:Chebyshev} 

As is well known, Dickson polynomials are closely related to Chebyshev polynomials.
The Chebyshev polynomials of the first kind, denoted $T_k(x)$, and Chebyshev polynomials of
the second kind, denoted $U_k(x)$,  belong to $\Z[x]$ and are given by the recursions:
\begin{eqnarray*} 
T_0(x)=1,&& T_1(x)=x,\qquad \,T_{k+2}(x) = 2x T_{k+1}(x) - T_k(x), \\
U_0(x)=1,&& U_1(x)=2x,\qquad \!\!U_{k+2}(x) = 2x U_{k+1}(x) - U_k(x). 
\end{eqnarray*}
These are related to Dickson polynomials of the first and second kind by the formulas
\begin{equation} D_k(2x)=2T_k(x),\qquad E_k(2x)=U_k(x). \label{DTEU} \end{equation}
Recall that the Dickson polynomials satisfy the functional equations
$$D_k(u+1/u) = u^k + 1/u^k,\qquad E_k(u+1/u) = \frac{u^{k+1}-1/u^{k+1}}{u-1/u},$$
where $k\ge 0$ and $u$ is an indeterminate. 
If we set $u=e^{i\theta}$, then $u+1/u=2\cos\theta$ and $u-1/u=2\sin\theta$.
Thus, the functional equations may be written as $D_k(2\cos\theta) = 2 \cos(k\theta)$ and $E_k(2\cos\theta) = \sin((k+1)\theta)/\sin\theta$.
This explains the famous functional equations for the Chebyshev polynomials that aid with trig substitutions 
in calculus:
$$T_k(\cos\theta) = \cos(k\theta),\qquad U_k(\cos\theta) \sin\theta = \sin((k+1)\theta).$$
From this, Schur observed that $T_k'(\cos\theta)(-\sin\theta)=-k\sin(k\theta)$,
therefore $T_k'(x)=kU_{k-1}(x)$.  In terms of Dickson polynomials, this is
equivalent to
\begin{equation} E_{k-1}(x)=(1/k) D_k'(x). \label{eq:DicksonDeriv} \end{equation}

Note that $D_k(x)$ and $E_k(x)$ are monic polynomials when $k\ge1$; whereas $T_k(x)$ has leading term $2^{k-1}x^k$
and $U_k(x)$ has leading term $2^kx^k$. Also, the Dickson polynomials have a characteristic~2 theory, while the Chebyshev
polynomials do not.   In that sense, the Dickson polynomials are more ``natural'' than the Chebyshev polynomials.
On the other hand, Chebyshev polynomials predate Dickson polynomials by about 40 years, they
have a wide range of applications, and thousands of articles have been written about them.
In this section, we restate theorems about Dickson polynomials in terms of their more famous cousin, the Chebyshev polynomials.

In odd characteristic, the formula (\ref{DTEU}) shows that $D_k$ and $T_k$ are conjugate: if $M_2$ denotes
the multiplication-by-2 map, then $T_k=M_2^{-1} \circ D_k \circ M_2$ and $U_k = E_k \circ M_2$.  
Thus, the theories of Dickson and Chebyshev polynomials are nearly identical in odd characteristic.

\begin{theorem} \label{thm:ChebyshevRefinements} ($q$ odd.) Let $\varep = \jacobi{-1}q$ and $\nu=\jacobi2q$, and let
$\calA_\lambda^{\varep_1,\varep_2}$ be as in (\ref{Aij}).
For any $k\ge 1$, $T_k(1)=1$, $T_k(-1)=(-1)^k$, and
\begin{enumerate}
\item[({\it i})] $T_k(\calA_1^{\nu,\nu}) \subset \calA_1^{\nu,\nu} \sqcup \{1,-\varep\}$ and 
$T_k(\calA_1^{-\nu,\nu}) \subset \calA_1^{-\nu,\nu}\sqcup \{1,\varep\}$.
\item[({\it ii})] If $k$ is odd then $T_k(\calA_1^{\varep\nu,-\nu}) \subset \calA_1^{\varep\nu,-\nu}\sqcup \{-1\}$
and $T_k(\calA_1^{-\varep\nu,-\nu}) \subset \calA_1^{-\varep\nu,-\nu}$.
\item[({\it iii})] If $k$ is even then $T_k(\calA_1^{-\nu,-\nu}) \subset \calA_1^{\nu,\nu}\sqcup \{1,-\varep\}$
and $T_k(\calA_1^{\nu,-\nu}) \subset \calA_1^{-\nu,\nu} \sqcup \{1,\varep\}$.
\item[({\it iv})] If $(q-1)/2$ divides $k$, then 
$$T_k(\F_q) \subset \calA_1^{+-} \sqcup \calA_1^{-+} \sqcup \{1,-1\}
= \left\{a \in \F_q : \text{$a^2-1$ is a nonsquare}\right\} \sqcup \{1,-1\}.$$ 
If in addition $k$ is even,
then $T_k(\F_q)\subset \calA_1^{-\nu,\nu} \sqcup \{1,-1\}$.
\item[({\it v})] If $(q+1)/2$ divides $k$, then 
$$T_k(\F_q)\subset \calA_1^{++} \sqcup \calA_1^{--} \sqcup \{1,-1\}
= \left\{a \in \F_q : \text{$a^2-1$ is a square}\right\}.$$ 
If in addition $k$ is even, then
$T_k(\F_q) \subset \calA_1^{\nu,\nu} \sqcup \{1,-1\}$.
\item[({\it vi})] If $\calA_1^{\varep_1,\varep_2}$ is nonempty, then $T_k$ permutes $\calA_1^{\varep_1,\varep_2}$ iff GCD$(k,d^{\nu\varep_1,\nu\varep_2})
= 1$, where $d^{++} = (q-1)/2$, $d^{+-}=q+1$, $d^{-+}=(q+1)/2$, $d^{--} = q-1$.
\end{enumerate}
\end{theorem}
\begin{proof}
Recall that $T_k(x)=(1/2)D_k(2x)$.  Then
$D_k(a)=a'$ iff $T_k(a/2) = a'/2$. Since $D_k$ respects the sets $\calA_2^{\varepsilon_1,\varepsilon_2}$,
we see that $T_k$ respects the sets $\left\{a/2 : a \in \calA_2^{\varepsilon_1,\varepsilon_2}\right\}$.
Note that
\begin{eqnarray*}
\left\{a/2 : a \in \calA_2^{\varepsilon_1,\varepsilon_2} \right\} &=& \left\{ b \in \F_q : \jacobi{2b-2}q=\varepsilon_1,\ 
\jacobi{2b+2}q= \varepsilon_2 \right\} \\
&=& \left\{ b \in \F_q : \jacobi{b-1}q=\varepsilon_1 \jacobi 2q,\  \jacobi{b+1}q=\varepsilon_2 \jacobi 2q \right\} \\
&=& \calA_1^{\nu\varep_1,\nu\varep_2}.
\end{eqnarray*}
Using these observations, Theorem~\ref{thm:ChebyshevRefinements} follows immediately from Theorem~\ref{thm:Aij}.
\end{proof}

\begin{theorem} ($q$ odd.) ({\it i}) If $(k,q^2-1)=1$, so that $D_k$ and $T_k$ permute $\F_q$, then the permutations of $\F_q$ that are induced by $D_k$ and $T_k$ 
are conjugate in the symmetric group on $\F_q$ and therefore have the same cycle structure. \\
({\it ii}) Let $\varep_1,\varep_2 \in \{1,-1\}$ and $\nu = \jacobi 2q$. If $(k,d^{\varep_1,\varep_2}) = 1$, so that $D_k$ permutes $\calA_2^{\varep_1,\varep_2}$ and
$T_k$ permutes $\calA_1^{\nu\varep_1,\nu\varep_2}$, then the cycle structures of these two permutations are the same.
\end{theorem}
\begin{proof} We observed in the paragraph preceding Theorem~\ref{thm:ChebyshevRefinements} that $T_k=M_2^{-1}\circ D_k \circ M_2$, where $M_2$ is the multiplication-by-2 map;
\ie, $T_k$ takes $a$ to $b$ iff $D_k$ takes $2a$ to $2b$. Thus, in the case that they are permutations,
a cycle $(a_1,a_2,\ldots,a_r)$ of $T_k$ corresponds to a cycle $(2a_1,2a_2,\ldots,2a_r)$ of $D_k$. This establishes ({\it i}). Noting that
$a\in\calA_1^{\nu\varep_1,\nu\varep_2} \iff 2a \in \calA_2^{\varep_1,\varep_2}$, the above observation also implies that ({\it ii}) holds.
\end{proof}

\begin{theorem} \label{thm:gijTU} ($q$ odd.)  Let $\varep=\jacobi{-1}q$, $\nu=\jacobi2q$, and $\calB_\lambda^{\varep_1,\varep_2}$ as in~(\ref{BDef}). 
Let $m=(q-\varep)/4$, which is an integer.  In $\F_q[x]$, the following holds. 
{} 

\medskip
\noindent
({\it i})\quad $T_m(x) = 2^{m-1} \prod\left\{x-b : b \in \calB_1^{-\nu,-\nu}\right\}$ and $U_{m-1}(x) = 2^{m-1}\prod \left\{ x-b : b \in \calB_1^{\nu,\nu} \right\}$.
({\it ii})\quad $U_{(q-3)/2}(x) = 2^{(q-3)/2} \prod \left\{ x - b : \jacobi{b^2-1}q = 1 \right \}$ and \\
$U_{(q-1)/2} = 2^{(q-1)2} \prod \left\{ x - b : \jacobi{b^2-1}q = -1 \right\}$.
\end{theorem}
\begin{proof} 
Since $T_k(x)=(1/2)D_k(2x)$ and $D_k$ is monic of degree~$k$, we see that $T_k$ has leading term $2^{k-1}x^k$, and if $D_k(x) = \prod_{b\in B} (x-b)$ then
$T_k(x) = 2^{k-1} \prod_{b \in B} (x-b/2)$.
Likewise, since $U_k(x)=E_k(2x)$ and $E_k$ is monic of degree~$k$, we see that $E_k$ has leading term $2^kx^k$ and if $E_k(x) = \prod_{b\in B'}(x-b)$ then
$U_k(x)=2^k \prod_{b\in B'}(x-b/2)$.
By Theorem~\ref{thm:fijDE}({\it iv}), $D_m(x) = \prod \{x-b : b \in \calB_2^{--}\}$ and $E_{m-1}(x) = \prod \{x-b : b \in \calB_2^{++} \}$.
Thus, $T_m(x)= 2^{m-1} \prod \{x-b : b \in (1/2) \calB_2^{--} \}$ and $U_{m-1}(x) = 2^{m-1} \prod \{ x-b : b \in (1/2)\calB_2^{++} \}$.
Since $(1/2)\calB_2^{\varep_1,\varep_2}=\calB_1^{\nu\varep_1,\nu\varep_2}$, ({\it i}) follows.
By Theorem~\ref{thm:fijDE}({\it iii}), $E_{(q-3)/2}(x) = \prod\{x-b : b \in \calS_q \}$ and $E_{(q-1)/2}(x) = \prod \{x-b : b \in \calN_q\}$,
where $\calS_q$ and $\calN_q$ are defined by (\ref{ZSN}). Then $U_{(q-3)/2}=2^{(q-3)/2} \prod\{x-b : b \in (1/2)\calS_q \}$ and
$U_{(q-1)/2} = 2^{(q-1)/2}\prod \{x-b : b \in (1/2) \calN_q \}$. Now 
$$(1/2)\calS_q = \left\{ a/2\in\F_q : \jacobi{a^2-4}q=1 \right\} = \left\{b\in\F_q : \jacobi{b^2-1}q=1 \right\}$$
and similarly $(1/2)\calN_q = \left\{b\in\F_q : \jacobi{b^2-1}q=-1\right\}$.  Part~({\it ii}) follows.
\end{proof}

\begin{proposition} If $q=p^n$ is odd, then $$U_{(q-1)/2}(x) + U_{(q-3)/2}(x) = \jacobi 2 q (x-1)^{(q-1)/2} \pmod p.$$
\end{proposition}
\begin{proof} By Corollary~\ref{cor:fSfNvalues2}, $E_{(q-1)/2}(x) + E_{(q-3)/2}(x) \equiv (x-2)^{(q-1)/2} \pmod p$.  Since $U_k(x)=E_k(2x)$,
\begin{equation*} U_{(q-1)/2}(x) + U_{(q-3)/2}(x) \equiv (2x-2)^{(q-1)/2} = \jacobi 2q (x-1)^{(q-1)/2} \pmod p.  \end{equation*}
\end{proof}

\section{Conclusions} \label{sec:conclusions}

A common theme in the study of finite fields is to describe a set
in multiple ways.  A simple example is $\{a \in \F_q^\x : a^k = 1 \} = \{ a^g : a \in \F_q^\x \}$, where $g = (q-1)/\GCD(q-1,k)$.
Two other examples are known as Hilbert's Theorem~90:
\begin{eqnarray}
\{ a \in \F_{q^e} : \Tr_{\F_{q^e}/\F_q}(a)=0 \} &=& \{ b - b^q : b \in \F_{q^e} \},\\
 \{a\in\F_{q^e}^\x : \Norm_{\F_{q^e}/\F_q}(a)=1 \} &=& \{ b^{q-1} : b \in \F_{q^e}^\x \}. \label{Hilbert90} 
\end{eqnarray}
Another example (due to Dillon and Dobbertin \cite{DD}) is $\delstar{q-1} = \calT_{0}$ and $\delstar{q+1} = \calT_{1}$ when $q=2^n$, where we recall $\del d = \{ a + 1/a : a^d = 1 \}$
and $\delstar d = \del d \setminus \{2,-2\}$ when $(d,q)=1$.
We note that if $q=2^n$ and $1\ne r\in \mu_{q-1}$ then
$$\frac 1 {r + 1/r} = \frac r {r^2+1} = \frac{ (r+1)+1 }{(r+1)^2} = \frac 1{r+1} + \frac 1{(r+1)^2};$$
this gives another way to see that $\Tr(1/y)=0$ when $0\ne y \in \del{q-1}$.

In similar vein, Brewer showed that
$\del{q-1} = \{2,-2\} \sqcup \calS_q$ and $\del{q+1} = \{2,-2\} \sqcup \calN_q$ when $q$ is odd, where
$\calS_q$ and $\calN_q$ are defined by (\ref{ZSN}).
We showed in Theorem~\ref{thm:AijDescription} that
$$\delstar{(q-1)/2} = \left\{ a \in \F_q : \jacobi{a-2}q = \jacobi{a+2}q = 1 \right\}$$
$$\delstar{(q+1)/2} = \left\{ a \in \F_q : \jacobi{a-2}q = -1, \jacobi{a+2}q = 1 \right\}$$
$$\delstar{(q+1)} \setminus \delstar{(q+1)/2}
= \left\{ a \in \F_q : \jacobi{a-2}q = 1, \jacobi{a+2}q = -1 \right\} $$
$$\delstar{(q-1)} \setminus \delstar{(q-1)/2}
= \left\{ a \in \F_q : \jacobi{a-2}q = \jacobi{a+2}q = -1 \right\}.$$ 
The above four  sets are denoted by $\calA_2^{++}$, $\calA_2^{-+}$, $\calA_2^{+-}$, and $\calA_2^{--}$, respectively.

Since the Dickson polynomial of the first kind satisfies $D_k(u+1/u)=u^k + 1/u^k$, it maps $\del d$ to $\del {d/g}$, where $g = \GCD(d,k)$. Because of
this, these polynomials respect the above sets in structured ways.  For example, $D_k$ always maps
$\{ a \in \F_q : \text{$a+2$ and $a-2$ are squares} \}$ to itself, and $D_k$ permutes that set iff GCD$(k,(q-1)/2)=1$. 
See Theorems~\ref{thm:oddq} and \ref{thm:Aij}. 
Our results about Dickson polynomials translate easily into new results about Chebyshev polynomials; see Section~\ref{sec:Chebyshev}.

Another result of this article is that
\begin{eqnarray*} D_m(x) &=& \prod\left \{ x - b : \jacobi{2-b}q=\jacobi{2+b}q = -1 \right\}\\
E_{m-1}(x) &=& \prod\left \{ x - b : \jacobi{2-b}q=\jacobi{2+b}q = 1 \right\},
\end{eqnarray*}
where $m= \lfloor (q+1)/4  \rfloor$.
By evaluating the first identity at $x=2$ and $x=-2$, we obtain
the ``Wilson-like theorems'' that were presented  in Section~\ref{sec:wilson-like}.
Many more formulas of this type are proved in the sequel to this article, \cite{wilson-like}.
Also, since closed-form expressions are known for the coefficients of $D_m(x)$ and $E_{m-1}(x)$, the above equations
yield explicit formulas for $\sigma_j^+$ and $\sigma_j^-$, where $\sigma_j^\mu$ is $j$th elementary symmetric function of
$\left\{ b \in \F_q : \jacobi{2-b}q=\jacobi{2+b}q = \mu \right\}$.
See Theorem~\ref{thm:sigmak}.

Finally, we found some curious identities such as (\ref{jacc}), (\ref{sqrtc}), and
Corollary~\ref{cor:fSfNvalues2}.

\section{Alternate proofs}  \label{sec:oneline}
Some results of this article are proved using Dickson polynomials, but may be stated without them.
One can ask whether there is an alternate proof that uses cleverness instead of Dickson polynomials.
Here are two such examples that were communicated to the author.

\subsection{Richard Stong.} 

Theorem~\ref{thm:wilson-like} uses Dickson polynomials to prove that for odd~$q$,
$$\prod\left\{a\in\F_q^\x : \text{$a$ and $4-a$ are nonsquares} \right\}=2.$$  
This alternative proof was proffered by Richard Stong:

Start with the product in question {\it except} that it runs over the residues:
\begin{equation*} Q = \prod\left\{ b\in\F_q^\x : \jacobi bq = \jacobi {4-b}q = 1 \right\}. \end{equation*}
Write $b=x^2=4-y^2$, where $x,y\in \F_q^\x$. Note that $x,y\ne0$ implies $x \not\in \{0,2,-2\}$
and $y \not\in\{0,2,-2\}$. Let $R$ denote a complete set of representatives for $(\F_q\setminus\{0,2,-2\})/\{\pm1\}$.
\begin{eqnarray*} Q &=& \prod\left\{ 4-y^2 : y\in R, \jacobi{4-y^2}q=1 \right\}\\
&=&  \prod\left\{ (2-y)(2+y) : y \in R, \jacobi{2-y}q=\jacobi{2+y}q \right\} \\
&=& \prod \left\{ 2-z : z \in \F_q \setminus \{0,2,-2\}, \jacobi{2-z}q=\jacobi{2+z}q \right\} \\
&=& \prod \left\{ c : c \in \F_q\setminus\{2,0,4\},\jacobi cq = \jacobi{4-c}q \right\} \\
&=& (1/2) \prod \left\{ c : c \in \F_q^\x,\ \jacobi cq = \jacobi{4-c}q \right\} \\
&=& (1/2) Q N,
\end{eqnarray*}
where $Q$ is the product over $c$ that are squares and $N$ is the product over $c$ that are nonsquares.
So $Q=(1/2)QN$, and $N=2$. Interestingly, this does not solve for $Q$.  It is shown
in \cite{wilson-like} that $Q=-\jacobi{-1}q/4$.

\subsection{John Dillon.}
John Dillon offered the following simple proof of the formula~(\ref{sqrtc}). Let $q=2^n$, and let $\calT_j=
\{a\in\F_q^\x: \Tr(1/a) = j \}$. Let $T_j(x) = \prod\{(x-a) : a \in \calT_j \}$. We are to show that $T_0(c)+T_1(c)=c^{1/2}$
for all $c\in\F_q$.  It suffices to show that $T_0(x)+T_1(x)=x^{q/2}$. Define
$$S_0(x)=\Tr(x) = x + x^2 + x^4 + \cdots x^{q/2},\qquad S_1(x) = S_0(x)+1.$$
Let $S_0^*(x)=S_0(x)/x$.
Note that $S_0^*(x)=0$ iff $x\in\F_q^\x$ and $\Tr(x)=0$, and $S_1(x)=0$ iff $x\in\F_q^\x$ and Tr$(x)=1$.  Thus,
$T_0(x)=\widetilde{S_0^*}(x)$ and $T_1(x)=\widetilde{S_1}(x)$, where $\widetilde f$ denotes the reverse of~$f$: 
$$\widetilde f(x) = x^{\deg(f)}f(1/x).$$ 
We have
\begin{eqnarray*}
T_0(x) + T_1(x) &=& x^{(q/2)-1} S_0^*(1/x) + x^{q/2} S_1(1/x) \\
&=& x^{(q/2)-1} S_0^*(1/x) + x^{q/2} (1 + S_0(1/x) ) \\
&=& x^{(q/2)-1} S_0^*(1/x) + x^{q/2} (1 + (1/x) S_0^*(1/x) ) \\
&=& x^{q/2}.
\end{eqnarray*}

\section{Acknowledgments}
The author thanks John Dillon and the anonymous reviewer for providing valuable comments that improved the article.

\end{document}